\documentclass[11pt, a4paper, reqno]{amsart}
\usepackage[utf8]{inputenc}
\usepackage[margin=1in,includehead,headheight=15pt,headsep=12pt]{geometry}
\usepackage{lmodern}
\usepackage{upgreek}
\usepackage{xcolor}
\usepackage{amsfonts}
\usepackage{amsmath}
\usepackage{amssymb}
\usepackage{amsthm}
\usepackage{amsrefs}
\usepackage{comment}
\usepackage{etaremune}
\usepackage{enumitem}
\usepackage{pdfpages}
\usepackage{fancyhdr}
\usepackage{setspace}
\usepackage{hyperref}
\usepackage{mathrsfs}
\usepackage[perpage]{footmisc}
\usepackage{hyperref}
\usepackage{url}
\usepackage{bbm}
\hypersetup{
	colorlinks=true,
	linkcolor=blue,
	filecolor=blue,      
	urlcolor=blue,
	citecolor=blue,
	pdftitle={},
	pdfpagemode=FullScreen,}

\numberwithin{equation}{section}

\newtheorem{thm}{Theorem}[section]
\newtheorem{lem}[thm]{Lemma}

\newtheorem{cor}[thm]{Corollary}

\theoremstyle{definition}
\newtheorem{re}[thm]{Remark}

\newcommand{\mb}[1]{\mathbb{#1}}
\newcommand{\st}[1]{\substack{#1}}
\newcommand{\lr}[1]{\left(#1\right)}

\newcommand{\p}[1]{\Phi_a{(#1)}}

\title{The distribution of $\alpha p$ modulo one in short intervals}

\author{Habibur Rahaman}

\address{Theoretical Statistics and Mathematics Unit, Indian Statistical Institute, 203 B.T. Road, Kolkata 700108, India}
\email{habibr0805@gmail.com}

\date{}

\begin{document}
\begin{abstract}
Let $\|x\|$ denote the distance from a real number $x$ to the nearest integer, and let $\alpha$ be any irrational number and $\beta_0$ any real number. We prove that there exist infinitely many $X$ such that for $Y=X^\gamma/2$ with $7/12<\gamma\le 1$ and $0<\theta<\min \{(12\gamma-7)/24, 1/9\}$, there is a prime $p$ satisfying
$$
\|\alpha p+\beta_0\|<X^{-\theta},\quad X-Y<p<X.
$$
Moreover, we obtain an asymptotic formula for the number of such primes. In particular, if $\alpha$ has a continued fraction expansion with bounded partial quotients, then this asymptotic holds for all sufficiently large $X$. Our result improves upon a previous result of Baier and Roy by extending the admissible range of the short interval parameter from $2/3<\gamma\le 1$ to $7/12<\gamma\le 1$, reaching Huxley's classical range for primes in short intervals, and also improves the exponent $\theta$ for $2/3\le\gamma\le3/4$. The improvement is obtained by adapting a recent result of Maynard, Pandey, and Radziwi{\l}{\l} on exponential sums over primes.
\end{abstract}
\maketitle
\section{Introduction}
Let $\alpha$ be any irrational number. By Dirichlet approximation's theorem (see, \cite[Theorem~1A]{Schmidt}), for any $Q\ge 2$, there exist coprime integers $a, q$ with $1\le q\le Q$ and 
\begin{align}\label{Dirichlet}
    \left|\alpha-\frac{a}{q}\right|\le\frac{1}{qQ}\le \frac{1}{q^2}.
\end{align}
This gives
\begin{align}\label{DA}
    \|\alpha q\|<q^{-\theta},
\end{align}
for infinitely many positive integers $q$ with $\theta=1$, where $\|x\|$ denotes the distance from $x$ to the nearest integer. There are many interesting results concerning Diophantine approximation when the denominators $q$ are restricted to some sparse subsets of the set of positive integers $\mb{N}$. For the most interesting case, i.e., when the denominators $q$ are restricted to primes, Vinogradov \cite{Vinagradov}, using his method of exponential sums over primes, first proved \eqref{DA} for infinitely many primes $p$ with $\theta<1/5$. This was later extended by several authors, including Vaughan \cite{Vaughan}, Harman \cite{Harman1,Harman2}, Heath-Brown and Jia \cite{HJ2002} and Matom\"aki \cite{Matomaki}, who proved $\theta<1/3$, which is the best known result. In fact, Matom\"aki's method gives for infinitely many $x$,
\begin{align*}
    \sum_{\st{x/2<p\le x\\\|\alpha p\|<x^{-\theta}}}1\gg \frac{x^{1-\theta}}{\log x},
\end{align*} for any $\theta<1/3$, which is of the expected order of magnitude.

On the other hand, a fundamental question in the theory of primes is to understand the distribution of primes in short intervals. If $\pi(x)$ denotes the prime counting function, then extending Montgomery's result Huxley \cite{Huxley} proved that $$\pi(x)-\pi(x-y)\sim \frac{y}{\log x}, \quad \text{ as } x\to\infty,$$
uniformly for $x^{7/12}<y\le x.$ Huxley's result remained the best until recently, when it was extended to $x^{17/30}<y\le x$ by Guth and Maynard in the groundbreaking paper \cite{GM}, where they proved new large value estimates for the Dirichlet polynomials and as an application of these estimates, they improved Huxley's exponent $7/12$ to $17/30$. It is also known that under the Riemann Hypothesis, the range of $y$ can be extended to $x^{1/2}<y\le x.$ If, instead of an asymptotic formula, one only seeks a lower bound of the correct order of magnitude, the best result is due to Baker, Harman and Pintz \cite{BHP}, who showed that $$\pi(x)-\pi(x-y)\gg \frac{y}{\log x},\quad \text{ as } x\to \infty,$$ uniformly for $x^{21/40}<y\le x$.

It is natural to consider the hybrid problem in which one asks for primes $p$ lying in a short interval and at the same time satisfying a Diophantine approximation condition. Recently, this problem was studied by Baier and Roy \cite{BR2026}, who proved that for a suitable function $F(t)$ (see the precise definition in \eqref{def_F} below) approximating the indicator function $1_{[0, x^{-\theta}]}(\|t\|)$, we have 
\begin{align}\label{BR1}
    \sum_{x-y<n\le x}\Lambda(n)F(n\alpha)\sim \delta y,\quad \text{as }x\to\infty,
\end{align}
for infinitely many $x$ depending on $\alpha$ and for $x^{2/3}<y\le x$, where $\delta=x^{-\theta}$ for some $\theta>0$, tending to zero as $y$ approaches $x^{2/3}.$ In particular, if $\alpha$ has a continued fraction expansion with bounded partial quotients, then \eqref{BR1} holds for all sufficiently large $x$, independently of $\alpha$. Moreover, as a corollary, they proved that 
\begin{align*}
    \sum_{\st{x-y<p< x\\\|p\alpha\|<\delta}}1\sim \frac{2\delta y}{\log x}\quad \text{as }x\to \infty,
\end{align*}
uniformly for $x^{2/3}<y\le x$, where the dependence of $x$ on $\alpha$ is as above. 
\subsection{Statement of the results}
The goal of this paper is to extend the hybrid result of Baier and Roy~\cite{BR2026} to the classical range of Huxley, $x^{7/12}<y< x$. To state our result, we first fix some parameters. By Dirichlet's approximation theorem, there are infinitely many $q$ satisfying \eqref{Dirichlet}. We fix such a denominator $q$ and define
\begin{align}\label{parameters}
    Y=X^{\gamma}/2,\ \gamma\le 1;\quad \delta=X^{-\theta}, \theta>0,
\end{align}
where 
\begin{align}\label{X,q relation}
    X\asymp q^{\frac{2}{2\gamma+\theta-1+\varepsilon}},
\end{align}
with $\varepsilon>0$ is sufficiently small and fixed.
We prove the following asymptotic formula which counts the number of primes $p\in[X-Y, X]$ with $X^{7/12}<Y\le X/2$, which also satisfy $\|\alpha p+\beta_0\|<X^{-\theta}$, for any fixed real number $\beta_0$.
\begin{thm}\label{thm_counting_primes}
    Let $\alpha$ be any irrational number and $\beta_0$ be any real number. Let $X, Y, \delta$ be defined as in \eqref{parameters} and \eqref{X,q relation}. Then, for infinitely many $X$, we have 
    \begin{align*}
        \sum_{\st{X-Y<p<X\\ \|\alpha p+\beta_0\|<\delta}}1\sim \frac{2\delta Y}{\log X}, \quad \text{ as } X\to \infty,
    \end{align*} provided 
    \begin{align*}
        \frac{7}{12}+2\theta+10\varepsilon<\gamma\le1,\quad \text{and}\quad\theta<\frac{1}{9}-\varepsilon.
    \end{align*}
\end{thm}
\begin{re}
For $\gamma=2/3,\ \beta_0=0$  and $\varepsilon$ arbitrarily small, our result gives $\theta<1/24$. In particular, this yields a non trivial approximation exponent $\theta$ already at $\gamma=2/3$, whereas the result of Baier and Roy~\cite{BR2026} applies in the range $\gamma>2/3$. More generally, a comparison of the admissible ranges of $\theta$ shows that our result is stronger for $\frac{7}{12}<\gamma<\frac{3}{4}$, while the result of Baier and Roy gives a larger admissible range of $\theta$ for $\frac{3}{4}<\gamma\le1$. Thus, by combining the two results, one could formulate a stronger piecewise statement covering the full range of $\gamma$. We do not pursue this formulation here, since our primary objective is to get shorter intervals, namely to extend the admissible range of $\gamma$, rather than to optimize the exponent $\theta$ over the full range of $\gamma$. 
\end{re}
\begin{cor}
Let $\alpha$ be any irrational number with bounded partial quotients in its continued fraction expansion and let $\beta_0$ be any real number. Let $X, Y, \delta$ be defined as in \eqref{parameters} and \eqref{X,q relation}. Then, for sufficiently large $X>0$, we have 
\begin{align*}
        \sum_{\st{X-Y<p<X\\ \|\alpha p+\beta_0\|<\delta}}1\sim \frac{2\delta Y}{\log X}, \quad \text{ as } X\to \infty,
    \end{align*} provided 
    \begin{align*}
        \frac{7}{12}+2\theta+10\varepsilon<\gamma\le1,\quad \text{and}\quad\theta<\frac{1}{9}-\varepsilon.
    \end{align*}
\end{cor}
\begin{proof}
Suppose that $\alpha$ is an irrational number with continued fraction expansion $\alpha=[a_0,a_1,\ldots]$, 
where the partial quotients are bounded, say $a_k\le M$ for all $k\ge 0$. Let $p_k/q_k$ denote the convergents of $\alpha$. Then $q_0=1,\ q_1=a_1$ and $q_k=a_kq_{k-1}+q_{k-2}$ for $ k\ge 2$. Consequently, $\{q_k\}_{k\in \mb{N}}$ is an increasing sequence of positive numbers such that $q_{k-1}<q_k\le (M+1)q_{k-1}$, for all $ k\ge 1$. Therefore, for sufficiently large $X>0$, we can always find some $q_n$ such that $q_n\asymp X^{(2\gamma+\theta-1+\varepsilon)/2}$, satisfying \eqref{X,q relation}. Hence the corollary follows.
\end{proof}
For technical reasons, we work with the following smoothed version and from that we deduce the above counting theorem.
For $x\in\mb{R}$, define
\begin{align}\label{def_F}
    F(x):=\sum_{m\in\mathbb Z}\exp\left(-\pi\frac{(x-m)^2}{\delta^2}\right).
\end{align}
Note that $F$ satisfies 
\begin{align*}
    \begin{cases}
        F(x)\gg 1, & \|x\|\in [0,\delta),\\
        F(x)\ll X^{-2026}, & \|x\|\in [X^\varepsilon\delta, 1/2],\\
        F(x)\ll 1, & x\in\mb{R},
    \end{cases}
\end{align*}
where $\varepsilon>0$ is arbitrarily small.

\noindent We consider the sum
\begin{align}\label{def_S}
    S:=\sum_{X-Y<n<X}\Lambda(n)F(n\alpha+\beta_0).
\end{align}
Our aim is to obtain an asymptotic formula for $S$. Again, to use analytic tools, it is more convenient to work with a suitable smooth approximations to the characteristic function of the interval $[0,1]$.

\noindent Let $A$ be a sufficiently large fixed constant and let 
\begin{align}\label{def_rho}
    \rho:=(\log X)^{-A}.
\end{align}
Also, let $w^-,w^+:\mb{R}\to[0,1]$ be smooth functions such that
\begin{align}\label{weights_order}
    w^-(u)\leq \mathbbm{1}_{[0,1]}(u)\leq w^+(u),
\end{align}
with
\begin{align}\label{weights_support}
    &w^-(u)=1 \quad\text{for }u\in[\rho,1-\rho],\qquad\text{supp}(w^-)\subseteq[0,1],\\
    &w^+(u)=1 \quad\text{for }u\in[0,1],\qquad\text{supp}(w^+)\subseteq[-\rho,1+\rho].
\end{align}
These functions can be chosen so that, for every fixed integer $j\ge 0$, the $j$-th derivative
\begin{align}\label{weights_derivatives}
    \lr{w^\pm}^{(j)}(u)\ll_j \rho^{-j} \ll_j (\log X)^{Aj}.
\end{align}
For $w\in\{w^+, w^-\}$, we define
\begin{align}\label{def_Sw}
    S(w):=\sum_{n\in\mb{N}}\Lambda(n)w\lr{\frac{X-n}{Y}}
    F(n\alpha+\beta_0).
\end{align}
Since \(\Lambda(n)F(n\alpha+\beta_0)\geq0\), it follows immediately from
\eqref{weights_order} that
\begin{align}\label{sandwich_S}
    S(w^-)\leq S\leq S(w^+).
\end{align}
We prove the following asymptotic formula for $S(w)$.
\begin{thm}\label{main_theorem}
Let $\alpha$ be any irrational number, and let $\varepsilon>0$ be arbitrarily small. Let $q,X,Y,\delta$ be as in \eqref{X,q relation} and \eqref{parameters}, and let $F$ be defined by \eqref{def_F}. Then, for $w\in\{w^+,w^-\}$ as above, we have
\begin{align}\label{smoothed_asymptotic}
S(w)=\delta Y+o(\delta Y),\quad X\to\infty,
\end{align}
for infinitely many \(X\), provided
\begin{align}\label{gamma1}
    \frac{7}{12}+2\theta+10\varepsilon<\gamma\le1,\quad \text{and}\quad \theta<\frac{1}{9}-\varepsilon.
\end{align}
\end{thm}

 In the work of Baier and Roy~\cite{BR2026}, Vaughan's identity is used to decompose the sums arising from the nonzero Fourier frequencies into Type I and Type II sums. With their choice of parameters, the Type I variable ranges up to $X^{2/3}$, while the Type II sums involve variables in the range from $X^{1/3}$ to $X^{2/3}$. The short interval range $Y>X^{2/3}$ is closely tied to these admissible Type I and Type II ranges. Our approach instead uses the recent ideas of Maynard, Pandey, and Radziwi{\l}{\l}~\cite{MPR} on exponential sums over primes, relating the nonzero Fourier frequencies to large-value estimates for products of Dirichlet polynomials. Suitably adapting these estimates to the short-interval setting allows us to extend the admissible range of $Y$ down to the exponent $7/12$.
\subsection{Notation} We will use the following notation throughout this paper.  
\begin{itemize}
    \item Expressions of the form ${f}(x)=O({g}(x))$, ${f}(x) \ll {g}(x)$, and ${g}(x) \gg {f}(x)$ signify that $|{f}(x)| \leq C|{g}(x)|$ for all sufficiently large $x$, where $C>0$ is an absolute constant. A subscript of the form $\ll_A$ means the implied constant may depend on the parameter $A$. Throughout the paper, the implied constants in $\ll$ and $\gg$ may depend at most on $\varepsilon, w$ and $F$, and such dependence will be suppressed most of the time. 
    \item The notation ${f}(x) \asymp {g}(x)$ indicates that ${f}(x) \ll {g}(x) \ll {f}(x)$.
    \item We will use $A\lessapprox B$, or $B\gtrapprox A$ to mean $A\ll X^{o(1)}B$, where $o(1)$ denotes a quantity that tends to zero as $X\to\infty.$
    \item  We denote $n\sim N$ to mean $N<n\leq 2N$.
    \item For any real number $x$, we will write $e(x):=e^{2\pi ix}$.
    \item We denote the Euler totient function by $\varphi$ and the divisor function by $d(n)$.
    \item Given any function $f\colon \mathbb{R}\to \mathbb{R}$, we define its Fourier transform $\widehat{f}\colon \mathbb{R}\to \mathbb{C}$ by
    \begin{align*}
        \widehat{f}(x):=\int_{-\infty}^\infty f(t)e(-xt)\: dt.
    \end{align*}
    \item We write $\chi(\bmod\, q)$ to denote a Dirichlet character modulo $q$.
    \item By $\Lambda(n)$ we denote the von Mangoldt function.
\end{itemize}
\subsection{Organization of the paper}
The paper is organized as follows. In Section~\ref{outline}, we outline the main idea of the proof of Theorem~\ref{main_theorem}. Section~\ref{preliminaries} contains the essential tools needed to prove the main theorem, including some elementary results on Fourier approximation and additive characters. In this section, we also obtain a short interval analogue of the Heath--Brown decomposition, collect results on Dirichlet polynomials, including Type II estimates, and carry out a stationary phase analysis to bound a Mellin transform. In Section~\ref{proof}, we use these tools to prove Theorem~\ref{main_theorem} and finally, in Section~\ref{Sec_count}, we deduce Theorem~\ref{thm_counting_primes} from Theorem~\ref{main_theorem}.
\subsection*{Acknowledgements} The author is thankful to Ritabrata Munshi for many helpful discussions and valuable suggestions. The author thanks Kunjakanan Nath for his encouragement in writing up this project and for his comments on an earlier draft. The author also thanks Soumya Bhattacharya and Stephan Baier for their comments on an earlier draft. The author is supported by the Research Associate Fellowship at the Indian Statistical Institute, Kolkata.

\section{Outline of the proof of Theorem~\ref{main_theorem}}\label{outline}
Let us sketch the proof of Theorem~\ref{main_theorem}. Recall that, for infinitely many $X$, we need to prove
$$
S(w)=\delta Y+o(\delta Y),\quad \text{as } X\to\infty,
$$
where $X,Y,\delta$ are given in \eqref{parameters} and satisfy \eqref{gamma1}. 

\noindent Using the Fourier expansion of $F(x)$ and truncating the Fourier series at $H=X^{o(1)}/\delta$, we can write, up to a negligible error (see Lemma~\ref{F_fourier} below) 
\begin{align*}
    S(w)= \delta\sum_{n\in\mb{
    N}}\Lambda(n)w\lr{\frac{X-n}{Y}}+\delta\sum_{0<|h|\leq H}e(h\beta_0)
    \exp\left(-\pi\delta^2h^2\right)\sum_{n\in\mb{N}}\Lambda(n)w\lr{\frac{X-n}{Y}}e(\alpha hn),
\end{align*}
where $H=X^{o(1)}/\delta$. Using the properties of the smooth function $w$ and the prime number theorem in short intervals due to Guth and Maynard (see Lemma~\ref{pnt2} below), the first term is 
$
\delta Y+o(\delta Y),
$
uniformly for $X^{17/30}<Y\le X$. (Though, for our range $X^{7/12}<Y< X$, Huxley's classical result is sufficient.) Therefore, to prove the main theorem, it is enough to show that the second sum, after removing the factor $\delta$, is $o(Y)$. 

We prove this by estimating the inner smoothed exponential sums over primes in short intervals
$$
S_h(w):=\sum_{n\in\mb{N}}\Lambda(n)w\lr{\frac{X-n}{Y}}e(\alpha hn),
$$
for each fixed $1\leq |h|\leq H$, and then summing trivially over $h$. To estimate these sums, we exploit the recent ideas of Maynard, Pandey, and Radziwi{\l}{\l} \cite{MPR}, who obtained a new bound for such exponential sums over primes without the short interval condition, improving the long-standing previous bound.

\noindent For any $R\ge 2$, Dirichlet's approximation theorem gives
\begin{align}\label{outline2}
    \alpha h=\frac{b_h}{r_h}+\eta_h
\end{align}
where $$(b_h, r_h)=1,\ r_h\le R\ \text{ and } \ |\eta_h|\le 1/r_hR.$$ 
Substituting this and using the representation of an additive character in terms of multiplicative characters (see Lemma~\ref{lem_of_character} below), we arrive roughly at the bound
\begin{align*}
    S_h(w)\lessapprox \frac{1}{r_h^{1/2}}\sum_{\chi (\bmod r_h)}\left|\sum_{n\in\mb{N}}\Lambda(n)\chi(n)e(n\eta_h)w\lr{\frac{X-n}{Y}}\right|.
\end{align*}
We define $$
G_h(u):=e(u\eta_h)w\lr{\frac{X-u}{Y}}.
$$
Then, using the Mellin inversion formula for $G_h(u)$, a smoothed version of \cite[Lemma~3.3]{MPR} that decomposes the Dirichlet series with coefficients $\Lambda(n)$ into a linear combination of suitable products of Dirichlet polynomials (see Lemma~\ref{lem_HB_identity} below), and the estimate for the Mellin transform $W_h(s)$ of $G_h(u)$ (see Lemma~\ref{lem_stationary} below), we roughly obtain
\begin{align}\label{outline1}
    S_h(w)\lessapprox\frac{Y}{Xr_h^{1/2}}\sup\limits_{j\ll (\log X)^{10}}\sum_{\chi (\bmod{r_h})}\int\limits_{-T_h}^{T_h}\left|F_j\lr{i(t+t_0), \chi}\right|\ dt,
\end{align}
where $t_0:=-2\pi \eta_h X,\ T_h= X^{\eta}\lr{\frac{X}{Y}+|\eta_h|Y}$ with $\eta=o(1)$, and each $F_j(s,\chi)$ are products of $k\le 10$ suitable Dirichlet polynomials. 

\noindent 
To estimate the inner average sum of Dirichlet polynomials, we use the recent results of Maynard, Pandey, and Radziwi{\l}{\l} (see Lemmas~\ref{lem_longfactor} and~\ref{lem_shortfactor} below). However, we must first verify that the parameters $r_h$ and $T_h$ satisfy the required size conditions. In particular, when $k\ge 5$, the application of Lemma~\ref{PNT_short} requires conditions of the form
$r_hT_h\gg X^{4/9}\ \text{ and }\ r_hT_h\lessapprox X^{1/2}$. We choose the parameter $R$ optimally so that $r_hT_h\lessapprox X^{1/2}$ for each $1\le |h|\le H$. To establish the complementary lower bound $r_hT_h\gtrapprox X^{4/9}$, we divide the set $\{h\in\mb{Z}: 1\le |h|\le H\}$ into two subsets according to the size of $r_h$. When $r_h$ is small, we compare the approximation \eqref{outline2} above with the fixed approximation 
$$\alpha=\frac{a}{q}+\beta,\ \text{  with}\quad (a,q)=1,\ |\beta|\le 1/q^2,$$
which gives a suitable lower bound for $r_h|\eta_h|$. When $r_h$ is comparatively large, we use its size directly. In this way, we obtain  $r_hT_h\gtrapprox X^{4/9}$ for all $1\le |h|\le H$, provided $\theta<1/9$. We can therefore apply Lemma~\ref{lem_combined} to bound the inner sum in \eqref{outline1}. This yields the required estimate provided $$\gamma>7/12+2\theta.$$  

\section{Preliminaries}\label{preliminaries}
In this section, we collect some of the tools needed to prove our main result. We begin with the following Fourier series expansion for $F$.
\begin{lem}\label{F_fourier}
Let $X, \delta$ be defined as in \eqref{parameters} and define
\begin{align}\label{def_H}
     H:= X^{\varepsilon}\delta^{-1}.
\end{align}
Then, uniformly for all $x\in\mathbb R$,
\[
    F(x)=\delta+\delta\sum_{1\le|h|\leq H}\exp\left(-\pi\delta^2h^2\right)e(hx)+O\left(X^{-2026}\right).
\]
\end{lem}
\begin{proof}
    See \cite[Proposition~5]{BR2026}.
\end{proof}
The following lemma is standard (for example, see \cite[Chapter~3]{IK}), which expresses the additive character in terms of multiplicative character.
\begin{lem}\label{lem_of_character}
    For $(n,q)=1$, we have
\begin{equation*}
    e\left(\frac{n}{q}\right)=\frac1{\varphi(q)}\sum_{\chi\,(\mathrm{mod}\, q)}\chi(n)\tau(\overline{\chi}),
\end{equation*}
where
\[
    \tau(\chi):=\sum_{b\, (\mathrm{mod}\,q)}\chi(b)e(b/q)
\]
is the Gauss sum modulo $q$ and $|\tau(\chi)|\le q^{1/2}$.  
\end{lem}
\begin{lem}[Prime number theorem in short intervals]\label{MG2026}
   Let $x,y\in\mb{R}$ and let $\varepsilon>0$ be arbitrarily small. As $x\to\infty$, we have 
    \begin{align*}
        \sum_{x-y<n<x}\Lambda(n)\sim y,
    \end{align*} for $x^{17/30+\varepsilon}<y<x/2$.
\end{lem}
\begin{proof}
   For $y\in [x^{17/30+\varepsilon}, x^{0.99}],$ this follows by an application of the recent result of Guth and Maynard on new large value estimate of Dirichlet polynomial (see \cite[Corollary~1.3]{GM}). For $y\in [x^{0.99}, x/2]$, it follows from Huxley \cite{Huxley}. 
\end{proof}
\begin{lem}\label{pnt2}
    let $\varepsilon>0$ be arbitrarily small and let $x,y\in\mb{R}$. Let $w\in\{w^-, w^+\}$, where $w^-, w^+$ are given in \eqref{weights_order} with $\rho=(\log x)^{-A},\ A\ge 2$. Then, we have 
    \begin{align*}
        \sum_{n\in\mb{N}}\Lambda(n)w\lr{\frac{x-n}{y}}=y+o(y),
    \end{align*}
    as $x\to\infty$, provided $x^{17/30+\varepsilon}<y\le x/2.$
\end{lem}
\begin{proof}
    First suppose $w=w^-$. Recall from \eqref{weights_support} that $w^-(u)=1,$ for $u\in[\rho, 1-\rho]$ and supported in $[0,1]$, where $\rho=(\log x)^{-A}$ given in \eqref{def_rho}.
Note that  
    \begin{align*}
        \sum_{n\in\mb{N}}\Lambda(n)w\lr{\frac{x-n}{y}}= \sum_{x-y\le n\le x}\Lambda(n)w\lr{\frac{x-n}{y}}
    \end{align*}
and 
\begin{align*}
   \sum_{x-(1-\rho)y\le n\le x-\rho y}\Lambda(n)\le  \sum_{x-y\le n\le x}\Lambda(n)w\lr{\frac{x-n}{y}}\le \sum_{x-y\le n\le x}\Lambda(n).
\end{align*}
Since 
\begin{align*}
    \sum_{x-(1-\rho)y\le n\le x-\rho y}\Lambda(n)&=\sum_{x-y\le n\le x}\Lambda(n)+O\lr{\sum_{x-y\le n\le x-(1-\rho)y}\Lambda(n)+\sum_{x-\rho y\le n\le x}\Lambda(n)}\\
    &=y+o(y)+O\lr{(\log x)\rho y}=y+o(y),
\end{align*}
the lemma follows for $w=w^-$. Similarly, the lemma also follows for $w=w^+$.
\end{proof}
\begin{lem}[Heath-Brown decomposition with a smooth weight]
\label{lem_HB_identity}
Let $x\ge 2$ be large and let $x^{4/9}\le B\le x$. Also, let $\chi$ be a Dirichlet character modulo $r$. Then for any smooth function $G$ compactly supported in $[x/2, 2x]$, and $c\in \mb{R}$, we have
\begin{align}\label{Left_HB}
    \sum_{n\in\mb{N}}\Lambda(n)\chi(n)G(n)
 =\frac{1}{2\pi i}\int\limits_{(c)}W(s)\lr{\sum_{j \le J}c_j  F_j(s,\chi)}ds,
 \quad J\ll (\log x)^{10},
\end{align}
where \begin{align*}
    W(s)=\int\limits_{0}^{\infty}G(u)u^{s-1} du,\quad s\in\mb{C},
\end{align*} the coefficients $|c_j|\ll1$ and
\begin{equation}
F_j(s,\chi)=\prod_{\ell=1}^{k_j}S_{\ell,j}(s, \chi),
\end{equation}
for some $k_j\le 10$ with 
\begin{align*}
    S_{\ell,j}(s,\chi)=\sum_{n\in I_{\ell,j}}a_{\ell,j}(n)\chi(n)n^{-s}.
\end{align*} 
Here $I_{\ell,j}\subseteq [N_{\ell,j}, 2^{10}N_{\ell,j}]$, for some $N_{\ell,j}\ge 1$ and coefficients $a_{\ell,j}(n)$ satisfying the following:
\begin{align*}
(1) &\ \prod_{\ell=1}^{k_j}N_{\ell,j}\asymp x,\\
(2) &\ N_{\ell_1,j}N_{\ell_2,j}\geq 2^{-20}B^{1/2}, \quad \text{ for }\ell_1\ne \ell_2,\\ 
(3) &\ a_{\ell,j}(n)\ll n^{o(1)}, \text{ for all } n\ge 1.\\
(4)& \text{ Moreover, if $N_{\ell,j}>B^{1/2}$, then $I_{\ell,j}=[N_{\ell,j},2N_{\ell,j}]$ and either the coefficients }
 a_{\ell,j}(n)=1\\& \text{ for all $n$ or }a_{\ell,j}(n)=\log n.
\end{align*}
\end{lem}
\begin{proof}
The proof follows along the same lines as that of \cite[Lemma~3.3]{MPR}. The main difference is that, instead of the sharp restriction $n\in[x,2x]$, we work with a smooth weight $G$ (which will allow us to localize to short intervals in our later application). This modification is handled by Mellin inversion and does not affect the Heath--Brown decomposition. For the sake of completeness, we provide some details. The Heath--Brown identity with parameter $K=5$ states that, for any $n\le2x$,
\begin{equation}
 \Lambda(n)=\sum_{j=1}^{5}(-1)^{j-1}\binom{5}{j}\sum_{\st{m_1m_2\cdots m_{2j}=n\\m_{j+1},\ldots,m_{2j}\le (2x)^{1/5}}}
 (\log m_1)\mu(m_{j+1})\cdots\mu(m_{2j}).
\end{equation}
Substituting this into the left-hand side of \eqref{Left_HB} gives
\begin{align*}
&\sum_{n\in\mb{N}}\Lambda(n)\chi(n)G(n)\\&
 =\sum_{j=1}^{5}(-1)^{j-1}\binom{5}{j}
 \sum_{\st{m_1,\ldots,m_{2j}\ge1\\
 m_{j+1},\ldots,m_{2j}\le(2x)^{1/5}}}
 (\log m_1)\mu(m_{j+1})\cdots\mu(m_{2j})\lr{\prod_{\ell=1}^{2j}\chi(m_\ell)}
 G(m_1\cdots m_{2j}).
\end{align*} 
We decompose each variable $m_\ell$ into dyadic ranges $N_\ell<m_\ell\leq2N_\ell,$ with $N_\ell$ a power of $2$, for $1\leq \ell\le 2j$ with $j\le 5$. Thus, there are at most ten variables and since $G$ is supported in $[x/2, 2x]$, the total number of resulting dyadic sums is $O((\log x)^{10})$. We may discard the part for which the product $m_1m_2\cdots m_{2j}$ does not lie in the support of $G$. On the remaining sums, we have
\[
 m_1\cdots m_{2j}\asymp x,
\]
and hence
\begin{align*}
    \prod_{\ell=1}^{2j}N_\ell\asymp x.
\end{align*}
Therefore, the left-hand side of \eqref{Left_HB} is equal to a linear combination of $O((\log x)^{10})$ terms of the form 
\begin{align*}
 \sum_{\st{m_1\sim N_1}}\sum_{\st{m_2\sim N_2}}\cdots\sum_{\st{m_{2j}\sim N_{2j}}}
 (\log m_1)\mu(m_{j+1})\cdots\mu(m_{2j})\lr{\prod_{\ell=1}^{2j}\chi(m_\ell)}
 G(m_1\cdots m_{2j}),
\end{align*}
where $N_\ell\le (2x)^{1/5}$ for $\ell= j+1, \ldots, 2j$ and 
\begin{align*}
    \prod_{\ell=1}^{2j}N_\ell\asymp x.
\end{align*}
We use the following  Mellin inversion formula 
\begin{align*}
    G(u)=\frac{1}{2\pi i}\int\limits_{(c)}W(s) u^{-s}\, ds , \ u>0.
\end{align*}
Note that, this is convergent for any $c\in \mb{R}$ by the fact that $G$ is smooth and compactly supported in $[x/2, 2x]$. Thus, the above sum reduces to
\begin{align*}
    \frac{1}{2\pi i}\int\limits_{(c)}W(s)\prod_{\ell=1}^{2j}S_\ell(s, \chi)\, ds,
\end{align*}
where $$S_\ell(s,\chi)=\sum_{m\sim N_\ell}\frac{a_\ell(m)\chi(m)}{m^s},$$
with $a_1(m)=\log m, \ a_\ell(m)=1,\ 2\le \ell \le j$ and $a_\ell(m)=\mu(m), \ j+1\le \ell\le 2j.$ The rest follows exactly in a similar way of the proof of \cite[Lemma~3.3]{MPR}, by merging two Dirichlet polynomials in each product $\prod_{\ell=1}^{2j}S_\ell(s, \chi)$, whose lengths satisfy $N_{\ell_1}N_{\ell_2}<2^{-20}B^{1/2}$ into a single Dirichlet polynomial and until no such pair remains. The coefficients of the resulting Dirichlet polynomials are finite Dirichlet convolutions of the original coefficients $1, \log n$, and $\mu(n)$, and hence satisfy the standard divisor bound $n^{o(1)}$. The remaining properties then follow exactly as in \cite[Lemma~3.3]{MPR}.
\end{proof}
The following key lemmas bound averages of products of Dirichlet polynomials.
\begin{lem}[Long factors]\label{lem_longfactor}
Let $x\geq 2$, $x^{4/9}\lessapprox B\leq x$, $1\leq T\ll x$, and
$1\leq q\leq x$ with $B\lessapprox qT$. Let $$F(s,\chi)=S_1(s,\chi)\cdots S_k(s,\chi),
\ k\leq 4,$$ be one of the products $F_j(s,\chi)$ appearing in Lemma~\ref{lem_HB_identity}. Then
\[
\sup_{\st{t_0\in\mb{R}\\|t_0|<x}}
\sum_{\chi\,(\bmod q)}
\int_{-T}^{T}
\left|F(i(t+t_0),\chi)\right|\,dt\lessapprox\lr{(qT)x^{1/2}+x}.
\]
\end{lem}
\begin{proof}
    This is \cite[Lemma~3.8]{MPR}.
\end{proof}
\begin{lem}[Short factors]\label{lem_shortfactor}
Let $x\geq 2$, $T\geq 1$, $q\geq 1$, and $B\gtrapprox x^{4/9}$ with $B\lessapprox qT\lessapprox x^{1/2}.$ Let $$F(s,\chi)=S_1(s,\chi)\cdots S_k(s,\chi),\ 5\le k\le 10,$$
be one of the products $F_j(s,\chi)$ appearing in Lemma~\ref{lem_HB_identity}. Then
\[
\sup_{t_0\in\mb{R}}\sum_{\chi\,(\bmod q)}\int_{-T}^{T}\left|F(i(t+t_0),\chi)\right|\,dt\lessapprox\lr{x+(qT)^{1/2}x^{19/24}}.
\]
\end{lem}
\begin{proof}
    This is \cite[Lemma~3.9]{MPR}.
\end{proof}
Combining the above two lemmas, we have the following lemma.
\begin{lem}\label{lem_combined}
Let $x\geq 2$, $T\geq 1$, $q\geq 1$, and $B\gtrapprox x^{4/9}$ with $B\lessapprox qT\lessapprox x^{1/2}$. Let $F(s,\chi)=S_1(s,\chi)\cdots S_k(s,\chi),
\ k\leq 10,$ be one of the products $F_j(s,\chi)$ appearing in Lemma~\ref{lem_HB_identity}. Then
\[
\sup_{\st{t_0\in\mb{R}\\|t_0|<x}}\sum_{\chi\,(\bmod q)}\int_{-T}^{T}
\left|F(i(t+t_0),\chi)\right|\,dt\lessapprox \lr{x+(qT)^{1/2}x^{19/24}}.
\]
\end{lem}
\begin{proof}
If $k\ge 5$, this is Lemma~\ref{lem_shortfactor}. Let $k\le 4$. Then, Lemma~\ref{lem_longfactor} gives
\begin{align*}
\sup_{|t_0|<x}\sum_{\chi\,(\bmod q)}\int_{-T}^{T}
\left|F(i(t+t_0),\chi)\right|\,dt\lessapprox \lr{x+(qT)x^{1/2}}.
\end{align*}
Since 
\begin{align*}
    (qT)x^{1/2}\lessapprox (qT)^{1/2}x^{3/4}< (qT)^{1/2}x^{19/24},
\end{align*}
the lemma follows.
\end{proof}
\subsection{Stationary phase analysis}
\begin{lem}\label{lem_stationary}
Let $q,X,Y,H$ be defined as in \eqref{parameters} and \eqref{def_H}, and let $ w\in\{w^+,w^-\} $ be one of the smooth weights defined as in \eqref{weights_order}. Let $\varepsilon>0$ be arbitrarily small. For each fixed $1\leq |h|\leq H$, define
$$
G_h(u):=w\lr{\frac{X-u}{Y}}e(u\eta_h),\qquad u>0,
$$ for some $\eta_h\in\mb{R}$ depending on $h$
and let
$$
W_h(s):=\int_0^\infty G_h(u)u^{s-1}\,du,\qquad s\in\mb{C},
$$
denote its Mellin transform. Set
\begin{align}\label{t0Th}
    t_0:=-2\pi\eta_hX,\qquad T_h:=X^\eta\lr{\frac{X}{Y}+|\eta_h|Y},
\end{align}
where $\eta=o(1)>0$ arbitrarily small and tending to zero as $X\to\infty.$
Then, for any fixed $C>0$,
$$W_h(it)\ll
    \begin{cases}
        \displaystyle
        \frac{Y}{X}\lr{1+\frac{|\eta_h|Y^2}{X}}^{-1/2},& |t-t_0|\ll T_h,\\
        X^{-C},& |t-t_0|\gg T_h.
    \end{cases}
$$
The implied constants may depend at most on $\eta, C$ and $A$.
\end{lem}
\begin{proof}
\textbf{Bound for small $t$}:
    Changing the variable as $v=(X-u)/Y$, we have 
    \begin{align}\label{W_h1}
        W_h(it)&=e(\eta_hX)Y\int\limits_{\mb{R}}w(v)e(-\eta_hYv)(X-vY)^{it-1} dv\nonumber\\
        &=e(\eta_hX)Y\int\limits_{\mb{R}}\frac{w(v)}{X-vY}e\lr{-\eta_hYv+\frac{t}{2\pi}\log(X-vY)} dv\nonumber\\
        &=e(\eta_hX)Y I(t),
    \end{align}
where 
\begin{align*}
    I(t):=\int_{\mb{R}}A(v)e\lr{B_t(v)} dv    
\end{align*}
with
\begin{align}\label{A and Bt}
        A(v):=\frac{w(v)}{X-vY},\qquad B_t(v):=-\eta_hYv+\frac{t}{2\pi}\log(X-vY).
    \end{align}
By integration by parts, 
\begin{align}\label{I(t)}
    I(t)\ll \sup\limits_{v_1\in[-\rho, 1+\rho]}|\Phi_t(v_1)|\int\limits_{\mb{R}}|A'(v)|dv\ll \frac{1}{X}\sup\limits_{v_1\in[-\rho, 1+\rho]}|\Phi_t(v_1)|,
\end{align}
where $$\Phi_t(v_1)=\int\limits_{-\rho}^{v_1}e(B_t(v)) dv.$$
Now $$B'_t(v)=-\eta_hY-\frac{tY}{2\pi (X-vY)}, \text{ and }\ B''_t(v)=\frac{-tY^2}{2\pi(X-vY)^2}.$$
So
\begin{align*}
    B'_t(v)=0\Longleftrightarrow t+2\pi\eta_h (X-vY)=0\Longleftrightarrow t=-2\pi \eta_h X+2\pi \eta_hYv.
\end{align*}
Since $-\rho\le v\le 1+\rho$, a stationary point can occur only if $|t-t_0|\leq 2\pi (1+\rho)|\eta_h|Y\ll |\eta_h|Y$ with $t_0=-2\pi \eta_hX$. 
Also, in this range, $$|B''_t(v)|=\frac{|t|Y^2}{2\pi(X-vY)^2}\asymp |\eta_h|Y^2/X.$$
Therefore, by the second derivative estimate (see \cite[Lemma~4.4]{Titchmarsh}), we have for $|t-t_0|\ll |\eta_h|Y$, 
\begin{align*}
    I(t)\ll \frac{1}{X\sqrt{|\eta_h|Y^2/X}}
\end{align*}
and together with \eqref{W_h1}, it gives 
\begin{align}\label{Wh2}
    W_h(it)\ll \frac{Y}{X}\lr{\frac{|\eta_h|Y^2}{X}}^{-1/2}.
\end{align}
\textbf{Bound for large $t$}:
Recall that,
\begin{align*}
    B'_t(v)=\frac{Y}{2\pi (X-vY)}\lr{-2\pi \eta_h X+2\pi \eta_hvY-t}=-\frac{Y}{2\pi(X-vY)}\lr{t-t_0-2\pi \eta_hvY}.
\end{align*}
Thus, 
\begin{align}\label{lower_B'}
    |B'_t(v)|\gg \frac{Y}{X}\lr{|t-t_0|-2\pi(1+\rho) |\eta_h|Y},
\end{align}
since $-\rho\le v\le 1+\rho$ and $Y\le X/2$.
As $\rho\le 1$, for $|t-t_0|\gg X^\eta \lr{\frac{X}{Y}+|\eta_h|Y}$, we have 
\begin{align*}
    |B'_t(v)|\gg X^\varepsilon.
\end{align*}
Define $D(v):=\frac{1}{B'_t(v)}$. 
Since
\[
\frac{d}{dv}e\lr{B_t(v)}
=2\pi iB'_t(v)e\left(B_t(v)\right),
\]
we have
\begin{align}\label{diff1}
e\left(B_t(v)\right)
=\frac{D(v)}{2\pi i}
\frac{d}{dv}e\left(B_t(v)\right).
\end{align}
Now, define recursively
\begin{align}\label{recurrence}
F_0(v):=A(v),\qquad F_{m+1}(v):=\frac{d}{dv}\lr{D(v)F_m(v)}, \quad m\ge0.
\end{align}
Applying \eqref{diff1} and integrating by parts once gives
\[
I(t)=-\frac{1}{2\pi i}\int_{-\rho}^{1+\rho}F_1(v)e\left(B_t(v)\right)\,dv.
\]
Because $w$ is compactly supported in $[-\rho,1+\rho]$, the boundary terms vanish. Repeating the same argument shows inductively that, for every integer $m\ge0$,
\begin{align}\label{multi_int_I}
I(t)=\frac{(-1)^m}{(2\pi i)^m}
\int_{-\rho}^{1+\rho} F_m(v)e\left(B_t(v)\right)\,dv.
\end{align}
We claim that, for every $m\ge1$,
\begin{align*}
F_m(v)=\sum_{r=0}^{m} \sum_{\substack{k_1,\ldots,k_m\ge0\\
r+k_1+\cdots+k_m=m}}c_m(r;k_1,\ldots,k_m)\, A^{(r)}(v)
\prod_{j=1}^{m}D^{(k_j)}(v),
\end{align*}
where the coefficients $c_m(r;k_1,\ldots,k_m)$ are non-negative integers depending only on $m$. We prove this claim in the next lemma.
Therefore, we have
\begin{align}\label{bound_Fm}
|F_m(v)|
\ll_m
\sum_{r=0}^{m}
|A^{(r)}(v)|
\sum_{\substack{r+k_1+\cdots+k_m=m\\k_j\ge0}}
\prod_{j=1}^{m}|D^{(k_j)}(v)|.
\end{align}
By Leibniz rule, \begin{align*}
    A(v)=\frac{w(v)}{X-vY}\Longrightarrow A^{(r)}(v)\ll_r \frac{\rho^{-r}}{X},
\end{align*}
and from
\begin{align*}
    D(v)=\frac{1}{B'_t(v)}=-\frac{2\pi(X-vY)}{Y}\cdot\frac{1}{t-t_0-2\pi \eta_h vY},
\end{align*}
we have for $r\ge 1$,
\begin{align*}
    D^{(r)}(v)
&=-\frac{2\pi}{Y}
\left
\{
(X-vY)
\lr{\frac{1}{t-t_0-2\pi \eta_h vY}}^{(r)}
-rY
\left(\frac{1}{t-t_0-2\pi \eta_h vY}\right)^{(r-1)}
\right\}\\
&=-\frac{2\pi}{Y}\left\{(X-vY)\frac{r!(2\pi \eta_hY)^r}{(t-t_0-2\pi \eta_hvY)^{r+1}}- rY \frac{(r-1)!(2\pi \eta_hY)^{r-1}}{(t-t_0-2\pi \eta_hvY)^r}\right\},
\end{align*}
which, for $|t-t_0|\gg X^{\eta}\lr{X/Y+|\eta_h|Y},$ gives
\begin{align*}
    |D^{(r)}(v)|\ll _r\frac{X}{Y}\cdot\frac{1}{|t-t_0|}, \quad r\ge 0.
\end{align*}
Therefore, from \eqref{multi_int_I} and \eqref{bound_Fm},
\begin{align*}
    I(t)\ll_m \int\limits_{-\rho}^{1+\rho}|F_m(v)|\, dv\ll_m \frac{1}{X}\lr{\frac{X\rho^{-1}}{Y|t-t_0|}}^m\ll_m \frac{1}{X}X^{-\eta m/2}\ll X^{-C-1},
\end{align*}
for any large $C$, choosing $m$ sufficiently large depending on $C, \eta$.
Therefore, for $|t-t_0|\gg X^\eta \lr{X/Y+|\eta_h|Y},$ we have from \eqref{W_h1}
$$W_h(it)\ll YX^{-C-1}\ll X^{-C},$$ for any large $C$.

\noindent\textbf{Bound for medium $t$}:
Let $|\eta_h|Y\ll |t-t_0|\ll X^{\eta}\lr{X/Y+|\eta_h|Y}$. Then from  the lower bound of first derivative \eqref{lower_B'},
\begin{align*}
   |B'_t(v)|\gg \frac{Y}{X}\lr{|t-t_0|-2\pi(1+\rho) |\eta_h|Y}\gg \frac{|\eta_h|Y^2}{X}.
\end{align*}
Thus, if $\frac{|\eta_h|Y^2}{X}\gg 1$, then by the first derivative bound for non-stationary points (see \cite[Lemma~4.2]{Titchmarsh}), we have from \eqref{W_h1} and \eqref{I(t)}, 
\begin{align}\label{Wh3}
    |W_h(it)|\ll \frac{Y}{X}\cdot\frac{X}{|\eta_h|Y^2}=\frac{1}{|\eta_h|Y}\ll \frac{1}{\sqrt{|\eta_h|X}}=\frac{Y}{X}\lr{\frac{|\eta_h|Y^2}{X}}^{-1/2}
\end{align}
and if, $|\eta_h|Y^2/X\ll 1$, then the trivial bound in \eqref{W_h1} gives 
\begin{align}\label{Wh4}
    |W_h(it)|\ll \frac{Y}{X}.
\end{align}
Thus, for $|t-t_0|\ll X^\eta \lr{X/Y+|\eta_h|Y}$,
we have from \eqref{Wh2}, \eqref{Wh3} and \eqref{Wh4},
\begin{align*}
    W_h(it)\ll \frac{Y}{X}\min\{1, (|\eta_h|Y^2/X)^{-1/2}\}\asymp \frac{Y}{X}\lr{1+\frac{|\eta_h|Y^2}{X}}^{-1/2}.
\end{align*} 
This completes the proof.
\end{proof}
\begin{lem}
    Let $F_m(v),\ m\ge 0$ be defined as in \eqref{recurrence}. Then for $m\ge 1$, we have 
\begin{align}\label{mth_derivative_F}
F_m(v)=\sum_{r=0}^{m} \sum_{\substack{k_1,\ldots,k_m\ge0\\
r+k_1+\cdots+k_m=m}}c_m(r;k_1,\ldots,k_m)\, A^{(r)}(v)
\prod_{j=1}^{m}D^{(k_j)}(v),
\end{align}
where the coefficients $c_m(r;k_1,\ldots,k_m)$ are non-negative integers depending only on $m$.
\end{lem}
\begin{proof}
We prove this lemma by induction on $m$. For $m=1$, \eqref{recurrence} gives
\[
F_1=A'D+AD',
\]
which is exactly of the form \eqref{mth_derivative_F}, the two terms correspond respectively to
\[
(r,k_1)=(1,0), \text{ and }(r,k_1)=(0,1).
\]
In both cases $r+k_1=1$.
Now assume that \eqref{mth_derivative_F} holds for some $m\ge1$.  Then a typical term of $F_m$ is of the form
\begin{align*}
T(v)=c_mA^{(r)}(v)
D^{(k_1)}(v)\cdots D^{(k_m)}(v),\quad \text{with }\ r+k_1+\cdots+k_m=m.
\end{align*}
By the recursion \eqref{recurrence}, this term contributes to $F_{m+1}$ through $\frac{d}{dv}\lr{D(v)T(v)}.$
Now, by Leibniz rule
\begin{align}\label{productrule}
\frac{d}{dv}\lr{D(v)T(v)}
=&c_mD'(v) A^{(r)}(v)\prod_{j=1}^{m}D^{(k_j)}(v)
+c_mD(v) A^{(r+1)}(v)\prod_{j=1}^{m}D^{(k_j)}(v)\nonumber\\
&\quad
+c_m\sum_{\ell=1}^{m}
D(v) A^{(r)}(v)
D^{(k_1)}(v)\cdots D^{(k_\ell+1)}(v)\cdots D^{(k_m)}(v).
\end{align}
Every term in \eqref{productrule} contains exactly $m+1$ factors of the form $D^{(k)}$.  Moreover, the total order of the derivative in each term is $m+1$. Indeed, in the first term the total order
of the derivative is
\[
1+r+k_1+\cdots+k_m=m+1,
\]
in the second term it is
\[
0+(r+1)+k_1+\cdots+k_m=m+1,
\]
and in the $\ell$-th term of the last sum in \eqref{productrule}, it is
\[
0+r+k_1+\cdots+(k_\ell+1)+\cdots+k_m=m+1.
\]
Thus every term in \eqref{productrule} has the form required in \eqref{mth_derivative_F} with $m$ replaced by $m+1$.  Combining equal monomials only changes their coefficients, which remain non-negative integers depending only on $m+1$. This proves \eqref{mth_derivative_F} by induction. 
\end{proof}
\section{Proof of Theorem~\ref{main_theorem}}\label{proof}
Substituting the Fourier expansion of $F$ from Lemma~\ref{F_fourier} into $S(w)$ and changing the order of summation in the second term, we obtain
\begin{align}\label{S(w)1}
S(w)
={}&
\delta\sum_{n\in\mb{N}}\Lambda(n)
w\lr{\frac{X-n}{Y}}
+\delta\sum_{0<|h|\leq H}e(h\beta_0)
\exp\left(-\pi\delta^2h^2\right)
\sum_{n\in\mb{N}}\Lambda(n)
w\lr{\frac{X-n}{Y}}e(\alpha hn)\nonumber\\&
+O(X^{-2025}).
\end{align}
By Lemma~\ref{pnt2}, we have
\begin{align}\label{PNT_short}
\delta\sum_{n\in\mb{N}}\Lambda(n)w\lr{\frac{X-n}{Y}}=\delta Y+o(\delta Y)=\delta Y(1+o(1)),
\quad \text{as } X\to\infty,
\end{align}
uniformly for $X^{17/30+\varepsilon}<Y\le X/2$. In particular, this holds throughout the range required in Theorem~\ref{main_theorem}. Thus, the theorem follows if we can show that the second term in \eqref{S(w)1} is $o(\delta Y)$.

Let $R\ge 2$. By Dirichlet's approximation theorem, there exist $b_h,\ r_h$ with $(b_h, r_h)=1$ and $1\le r_h\le R$ such that 
\begin{align}\label{Dirichlet3}
    \alpha h=\frac{b_h}{r_h}+\eta_h,
\end{align} 
with $|\eta_h|\le \frac{1}{r_h R}.$ For each fixed $h$ with $1\le |h|\le H$, we define 
\begin{align}\label{def_Sh(w)}
    S_h(w):=\sum_{n\in\mb{N}}\Lambda(n)w\lr{\frac{X-n}{Y}}e(\alpha hn).
\end{align}
Substituting the above value of $\alpha h$ from \eqref{Dirichlet3}, we can write
\begin{align*}
S_h(w)&=\sum_{n\in\mb{N}}\Lambda(n)w\lr{\frac{X-n}{Y}}e(nb_h/r_h)e(n\eta_h)\\
&=\sum_{(n,r_h)=1}\Lambda(n)w\lr{\frac{X-n}{Y}}e(nb_h/r_h)e(n\eta_h)+O(X^\varepsilon),
\end{align*}
as the total contribution of the terms with $(n,r_h)>1$ is bounded by 
\begin{align*}
    \sum_{\st{p|r_h}}\sum_{\st{j\\p^j\ll X}}\log p\ll X^\varepsilon.
\end{align*}
For the remaining $n$, as $(nb_h, r_h)=1$, it follows from Lemma~\ref{lem_of_character} that 
\begin{align*}
    S_h(w)=\frac{1}{\varphi(r_h)}\sum_{\chi (\bmod r_h)}\chi(b_h)\tau(\bar{\chi})\sum_{n\in\mb{N}}\Lambda(n)\chi(n)e(n\eta_h)w\lr{\frac{X-n}{Y}}+O(X^\varepsilon).
\end{align*}
Note that, we have removed the condition $(n, r_h)=1$ from the inner sum as it is guaranteed by $\chi(n).$
Define \begin{align*}
    G_h(u)=e(u\eta_h)w\lr{\frac{X-u}{Y}}, \ \text{ for }u\in\mb{R},
\end{align*}
with $\eta_h$ defined by \eqref{Dirichlet3}. Using Heath-Brown decomposition, Lemma~\ref{lem_HB_identity} to the inner sum with $x=X$ and using the bound of the Mellin transform $W_h(s)$ of $G_h(u)$ from Lemma~\ref{lem_stationary}, we have 
\begin{align}\label{bouns_sh}
    S_h(w)\ll X^{\varepsilon}\frac{r_h^{1/2}}{\varphi(r_h)}\cdot\frac{Y}{X}\lr{1+\frac{|\eta_h|Y^2}{X}}^{-1/2}\sup\limits_{j\ll (\log X)^{10}}\sum_{\chi (\bmod{r_h})}\int\limits_{-T_h}^{T_h}\left|F_j\lr{i(t+t_0), \chi}\right|\ dt,
\end{align}
where $t_0,\ T_h$ are given in \eqref{t0Th} and 
\begin{equation*}
F_j(s,\chi)=\prod_{\ell=1}^{k_j}S_{\ell,j}(s, \chi),
\end{equation*} 
are products of Dirichlet polynomials $S_{\ell,j}$ given in Lemma~\ref{lem_HB_identity}.

We want to use the estimate from Lemma~\ref{lem_combined} to bound the above $S_h(w)$. For that, we need to show the conditions $r_hT_h\lessapprox X^{1/2}$ and $r_hT_h\gtrapprox X^{4/9}$ hold, for all $1\le |h|\le H$.
Since, $r_h\le R$ and $|\eta_h|r_h\le 1/R$, we have
\begin{align*}
    r_hT_h=r_hX^\eta\lr{\frac{X}{Y}+|\eta_h|Y}\le X^\eta \lr{\frac{XR}{Y}+\frac{Y}{R}}.
\end{align*}
For optimality, we take $R=Y/X^{1/2}$, and hence 
\begin{align*}
    r_hT_h\ll X^{\frac{1}{2}+\eta}\lessapprox X^{1/2}.
\end{align*}
\subsection{Lower bound of $r_hT_h$}
To show $r_hT_h\gtrapprox X^{4/9}$, we next distinguish two ranges according to the size of the denominators $r_h$. The reason for doing so is that when $r_h$ is small, we show that the product $r_h|\eta_h|$ still has some size. On the other hand, when $r_h$ is large, we directly use the size of $r_h$. 

\noindent \textbf{Case~1:} Let $r_h\le q/2H$. By \eqref{Dirichlet}, we have
\[
\alpha=\frac{a}{q}+\beta,\qquad |\beta|\le q^{-2},
\]
that is, for $1\le h\le H$,
\[
h\alpha=\frac{a_h}{q_h}+h\beta,\qquad q_h=\frac{q}{(q,h)}\ge \frac{q}{H},
\]
where $a_h/q_h$ is the reduced form of $ha/q$. Thus $h\alpha$ already has a rational approximation whose denominator is at least $q/H$.
We can therefore compare the new approximation $b_h/r_h$ of $\alpha h$ with this inherited approximation $a_h/q_h$. 

\noindent Since $r_h\le \frac{q}{2H}<q/H\le q_h$, the two reduced fractions $b_h/r_h$ and $a_h/q_h$ are distinct and therefore
\[
\left|\frac{b_h}{r_h}-\frac{a_h}{q_h}\right|= \left|\frac{b_hq_h-a_hr_h}{r_hq_h}\right|\ge\frac1{r_hq_h}.
\]
Since
\[
\frac{b_h}{r_h}+\eta_h=\frac{a_h}{q_h}+h\beta,
\]
it follows that
\[
|\eta_h|\ge\left|\frac{b_h}{r_h}-\frac{a_h}{q_h}\right|-|h\beta| \ge\frac1{r_hq_h}-|h\beta|.
\]
Multiplying by $r_h$, and using $q_h\le q,\ |h|\le H$, and $r_h\le q/(2H)$, we obtain
\begin{align}\label{lower_bound_rn}
r_h|\eta_h|\ge\frac1{q_h}-\frac{r_hH}{q^2}\ge\frac1q-\frac1{2q}=
\frac1{2q}.
\end{align}
Thus, in the range $r_h\le \frac{q}{2H},$ the product $r_h|\eta_h|$ cannot be arbitrarily small. Also 
\begin{align*}
    r_hT_h=r_hX^\eta\lr{\frac{X}{Y}+|\eta_h|Y}\ge X^\eta r_h|\eta_h| Y\ge X^\eta Y/2q\gtrapprox X^{4/9},
\end{align*}
provided 
\begin{align}\label{con_q1}
    q\ll YX^{-\frac{4}{9}}=X^{\gamma-\frac{4}{9}}.
\end{align}
\noindent\textbf{Case~2:} Let $r_h>q/2H$. Then 
\begin{align*}
    r_hT_h=r_hX^\eta\lr{\frac{X}{Y}+|\eta_h|Y}\ge X^\eta Xq/2YH\gtrapprox X^{4/9},
\end{align*}
provided 
\begin{align}\label{con_q2}
    q\gg YHX^{\frac{4}{9}-1}\gg X^{\gamma+\theta-\frac{5}{9}+\varepsilon}.
\end{align}
From \eqref{con_q1} and \eqref{con_q2}, for the existence of such $q$, we should have $$\gamma+\theta-\frac{5}{9}+\varepsilon<\gamma-\frac{4}{9},$$ which is true for $\theta<1/9-\varepsilon$ and taking average of these powers, we see that our choice (see \eqref{X,q relation}) 
$$q\asymp X^{\gamma+\frac{\theta}{2}-\frac{1}{2}+\frac{\varepsilon}{2}},$$
satisfies both \eqref{con_q1} and \eqref{con_q2} for $\theta<\frac{1}{9}-\varepsilon$ and we have $r_hT_h\gtrapprox X^{4/9}$ for all $1\le |h|\le H.$
\subsection{The final estimate}
\noindent Invoking the bound from Lemma~\ref{lem_combined} with  $B=X^{4/9}$, and $X, r_h, T_h$ in place of $x, q, T$, respectively, and choosing $\eta>0$ suitably small depending on $\varepsilon$, we have from \eqref{bouns_sh},
\begin{align*}
    S_h(w)\ll X^{2\varepsilon}\frac{Y}{r_h^{1/2}X}\lr{1+\frac{|\eta_h|Y^2}{X}}^{-1/2}\lr{X+(r_hT_h)^{1/2}X^{19/24}}.
\end{align*}
Since 
$$T_h=X^\eta\lr{\frac{X}{Y}+|\eta_h|Y}=\frac{X^{1+\eta}}{Y}\lr{1+\frac{|\eta_h|Y^2}{X}},$$
this gives 
\begin{align}\label{Sh1}
    S_h(w)&\ll X^{2\varepsilon}\lr{\frac{Y}{r_h^{1/2}}\lr{1+\frac{|\eta_h|Y^2}{X}}^{-1/2}+Y^{1/2}X^{7/24+\eta/2}}.
\end{align}
For $r_h\le q/2H$, using $\lr{1+\frac{|\eta_h|Y^2}{X}}^{-1/2}\le\lr{\frac{|\eta_h|Y^2}{X}}^{-1/2}$ and \eqref{lower_bound_rn}, we have from \eqref{Sh1} 
\begin{align}\label{bound_sh_small_r}
    S_h(w)\ll X^{2\varepsilon}\lr{\frac{X^{1/2}}{(r_h|\eta_h|)^{1/2}}+Y^{1/2}X^{7/24+\eta/2}}
    \ll X^{2\varepsilon}\lr{X^{1/2}q^{1/2}+Y^{1/2}X^{7/24+\eta/2}}.
\end{align}
On the other hand, when $r_h>\frac{q}{2H}$, using $\lr{1+\frac{|\eta_h|Y^2}{X}}^{-1/2}\le 1$, we have from \eqref{Sh1}
\begin{align}\label{bound_Sh_large_r}
    S_h(w)\ll X^{2\varepsilon}\lr{\frac{Y}{r_h^{1/2}}+Y^{1/2}X^{7/24+\eta/2}}\ll X^{2\varepsilon}\lr{\frac{YH^{1/2}}{q^{1/2}}+Y^{1/2}X^{7/24+\eta/2}}.
\end{align}
Combining \eqref{S(w)1}, \eqref{PNT_short} and \eqref{def_Sh(w)}, we can write 
\begin{align*}
    S(w)=\delta Y+o(\delta Y)+\delta\hspace*{-5pt}\sum_{\st{1\le |h|\le H\\r_h\le q/2H}}e(h\beta_0)\exp(-\pi \delta^2 h^2)S_h(w)+\delta \hspace*{-5pt}\sum_{\st{1\le |h|\le H\\r_h>q/2H}}e(h\beta_0)\exp(-\pi \delta^2 h^2)S_h(w).
\end{align*}
Substituting the bounds from \eqref{bound_sh_small_r} and \eqref{bound_Sh_large_r} and summing over $|h|\le H=X^{\varepsilon}/\delta$ trivially, this implies 
\begin{align*}
    S(w)=\delta Y+o(\delta Y)+O\lr{X^{4\varepsilon}\lr{X^{1/2}q^{1/2}+Y^{1/2}X^{7/24}+YH^{1/2}/q^{1/2}}},
\end{align*}
where $\eta$ is chosen sufficiently small depending on $\varepsilon$.
Recall the choices of $X, Y, \delta$ from \eqref{X,q relation}, \eqref{parameters}. Then, writing every other parameters in terms of $X$, the first term in side the bracket 
\begin{align}\label{gma1}
    X^{4\varepsilon}X^{1/2}q^{1/2}\ll \delta Y/X^\varepsilon\Longleftrightarrow \gamma>\frac{5}{2}\theta+\frac{1}{2}+\frac{21}{2}\varepsilon,
\end{align}
the second term 
\begin{align}\label{gma2}
    X^{4\varepsilon}Y^{1/2}X^{7/24}\ll \delta Y/X^\varepsilon\Longleftrightarrow \gamma>\frac{7}{12}+2\theta +10\varepsilon,
\end{align}
and the third term
\begin{align}\label{gma3}
    X^{4\varepsilon}YH^{1/2}/q^{1/2}\ll \delta Y/X^\varepsilon\Longleftrightarrow \gamma>\frac{5}{2}\theta+\frac{1}{2}+\frac{21}{2}\varepsilon.
\end{align}
Also taking $\varepsilon$ sufficiently small, we see that the condition \eqref{gma2} dominates the other two conditions \eqref{gma1} and \eqref{gma3} for $\theta<1/9-\varepsilon$ and therefore, 
\begin{align*}
    S(w)=\delta Y +o(\delta Y),\quad \text{as }X\to\infty,
\end{align*}
provided 
\begin{align*}
    \gamma>\frac{7}{12}+2\theta +10\varepsilon,\quad \theta<\frac{1}{9}-\varepsilon.
\end{align*}
Since Dirichlet's theorem \eqref{Dirichlet} provides infinitely many denominators $q$, we may choose infinitely many corresponding values of $X$ satisfying \eqref{X,q relation}. Along this sequence, we can take $X\to\infty$, and hence the above argument makes sense. This completes the proof of Theorem~\ref{main_theorem}.

\section{Proof of Theorem~\ref{thm_counting_primes}}\label{Sec_count}
In this section, we prove Theorem~\ref{thm_counting_primes} from Theorem~\ref{main_theorem}. By Theorem~\ref{main_theorem}, we have for $w\in\{w^-, w^+\}$,
\begin{align*}
    S(w^\pm)=\delta Y\left(1+O(\rho)\right)+o(\delta Y)=\delta Y(1+o(1)).
\end{align*}
Together with \eqref{sandwich_S}, this gives
\begin{align*}
    \delta Y(1-o(1))\leq S\leq\delta Y(1+o(1)),
\end{align*}
that is,
$$S\sim \delta Y,\quad \text{as  } X\to\infty,$$
where $S$ is given in \eqref{def_S}. Also, in $S$, the total contribution of prime powers $p^k,\ k\geq2$, is $$O(X^{1/2}\log X)=o(\delta Y),$$
provided $$\gamma>\frac{1}{2}+\theta+\varepsilon,$$
which is true as $\gamma>\frac{7}{12}+2\theta+10\varepsilon.$ Therefore,
$$
\sum_{X-Y<p<X}(\log p)F(p\alpha)
\sim\delta Y.
$$
But $\log p\sim\log X$  uniformly for $X-Y<p<X$, and consequently
$$
\sum_{X-Y<p<X}F(p\alpha)\sim\frac{\delta Y}{\log X}.
$$
Finally, as pointed out by Baier and Roy~\cite{BR2026}, to pass from the smooth function $F$ to the sharp condition $\|\alpha p\|<\delta$, we may work with the standard Fourier majorant and minorant for the indicator function of $\|\alpha p\|<\delta$ (see \cite[Lemma~2.1]{PDS}). The zero frequency term will give the expected main term $2\delta Y/\log X$, while the resulting nonzero Fourier frequencies will lie in the same range as those occurring in the proof of Theorem~\ref{main_theorem} and can be absorbed into the error term. We may therefore omit the rest of the proof.
\bibliography{biblio}

\end{document}